\documentclass[10pt]{article}
\usepackage{amsmath}
\usepackage{amsthm}
\usepackage{amssymb}
\usepackage{mathptmx}
\usepackage{concrete}
\def\pf{{\bf Proof. }}
\def\df#1#2{\frac{\displaystyle #1}{\displaystyle #2}}
\def\tpds#1#2{\df{\partial #1}{\partial #2}}
\def\bel#1{\begin{equation}\label #1}

\def\U{\Upsilon}

\def\ro{\rho}

\def\nd{\noindent}

\def\U{\Upsilon}

\def\BC{\Bbb C}
\def\BCB{\overline{\Bbb C}}

\def\P2{{\Bbb C}P^2}

\def\SU{{\rm PSU}(2)}
\def\PGL{{\rm PGL}(2,\,{\Bbb C})}
\def\U{{\rm U}(1)}
\def\hf{\hfill{$\Box$}}
\def\pf{{\bf Proof.} \ }
\def\<{\leq}
\def\>{\geq}

\newtheorem{thm}{Theorem}[section]
\newtheorem{lem}{Lemma}[section]
\newtheorem{prop}{Proposition}[section]
\newtheorem{cor}{Corollary}[section]

\newtheorem{exam}{Example}[section]

\begin{document}

\title{\bf
Conformal metrics with constant curvature one and finite
conical singularities on compact Riemann surfaces}
\author{Qing Chen$^\dag$, Wei Wang, Yingyi Wu$^\ddag$ and Bin Xu$^\S$}
\maketitle

\nd {\small {\bf Abstract:}
A conformal metric $g$ with constant
curvature one and finite conical singularities on a compact
Riemann surface $\Sigma$ can be thought of as the pullback of the
standard metric on the 2-sphere by a multi-valued locally
univalent meromorphic function $f$ on $\Sigma\backslash \{{\rm
singularities}\}$, called the {\it developing map} of the metric
$g$.  When the developing map $f$ of such a metric $g$ on the
compact Riemann surface $\Sigma$ has reducible monodromy, we show
that, up to some M{\" o}bius transformation on $f$, the
logarithmic differential $d\,(\log\, f)$ of $f$ turns out to be an
abelian differential of 3rd kind on $\Sigma$, which satisfies some
properties and is called a {\it character 1-form of} $g$.
Conversely, given such an abelian differential $\omega$ of 3rd
kind satisfying the above properties, we prove that there exists a
unique conformal metric $g$ on $\Sigma$ with constant curvature
one and conical singularities such that one of its character
1-forms coincides with $\omega$. This provides new examples of
conformal metrics on compact Riemann surfaces of constant
curvature one and with singularities. Moreover, we prove that the
developing map is a rational function for a conformal metric $g$
with constant curvature one and finite conical singularities with
angles in $2\pi\,{\Bbb Z}_{>1}$ on the two-sphere.}     \\

\nd{\bf Keywords.} conformal metric of constant curvature one,
conical singularity, developing map, character 1-form\\

\nd{\bf 2010 Mathematics Subject Classification.} Primary 32Q30;
secondary 34M35

\footnotetext[1]{$^\dag$The first author is supported in part by
the National Natural Science Foundation of China (Grant No.
11271343).}

\footnotetext[2]{$^\ddag$The third author is supported in part by
the National Natural Science Foundation of China (Grant No.
11071249) and the President Fund of UCAS.}

\footnotetext[3]{$^\S$The last author is supported in part by
Anhui Provincial Natural Science Foundation (Grant No.
1208085MA01) and the Fundamental Research Funds for the Central
Universities (Grant No. WK0010000020).}

\section{Introduction}

Let $\Sigma$ be a compact Riemann surface and $p$ a point on
$\Sigma$. A conformal metric $g$ on $\Sigma$ has a {\it conical
singularity} at $p$ with the {\it singular angle} $2\pi\alpha>0$
if in a neighborhood of $p$, $g=e^{2\varphi}\,|dz|^2$, where $z$
is a local complex coordinate defined in the neighborhood of $p$
with $z(p)=0$ and $\varphi-(\alpha-1)\,\ln\,|z|$ is continuous in
the neighborhood. Let $p_1,\cdots,p_n$ be points on $\Sigma$ and
$g$ a conformal metric on $\Sigma$ with conical singularity at
$p_j$ of singular angle $2\pi\alpha_j>0$ for $j=1,\cdots,n$. Then
we call that the metric $g$ represents the divisor $
D:=\sum_{j=1}^n\,(\alpha_j-1)\, P_j$. The Gauss-Bonnet formula
says that the integral of the curvature on $\Sigma$ equals $2\pi$
times
\[\chi(\Sigma)+\deg\,D,\]
where $\chi(\Sigma)$ denotes the Euler number of $\Sigma$, and
$\deg\,D=\sum_j\,(\alpha_j-1)$ the degree of the divisor $D$. A
classical problem is whether there exists a conformal metric on
$\Sigma$ of constant curvature $K$ representing the divisor $D$.
If $K\leq 0$, then the unique metric exists if and only if the
left hand side $\chi(\Sigma)+\deg\,D\leq 0$; see \cite{MO88,
Tr91}.

If $\chi(\Sigma)+\deg\,D>0$, or equivalently $K\equiv 1$ if we
multiply the original metric by some constant, the problem turns
to be quite subtle and is still open now, except that there are
some partial results. Troyanov \cite{Tr89} considered the case of
two points on the sphere and proved that the necessary and
sufficient condition in this case is $\alpha_1=\alpha_2$. A more
general result also due to him \cite[Theorem 4]{Tr91} says that
there exists a metric of constant positive curvature if
\begin{equation}
\label{equ:tro}
0<\chi(\Sigma)+\deg\,D<\min\,\{2,\,2\min\,\alpha_j\}.
\end{equation}
Luo and Tian \cite{LT92} proved that the above condition is also
necessary and the metric is unique, provided that $\Sigma$ is the
2-sphere and all angles lie in $(0,\,2\pi)$. In case that $\Sigma$
is a sphere and the divisor $D$ supports at three points,
Umehara-Yamada \cite{UY00}, Eremenko \cite{Er04}, Furuta-Hattori
\cite{FH} and Fujimori-Kawakami-Kokubo-Rossman-Umehara-Yamada
\cite{FKKRUY11} give a necessary and sufficient condition for the
existence of the metric, which is also unique if and only if none
of the three angles belongs to
$2\pi\,{\Bbb Z}_{>0}$. \\

We attack the problem by using the idea of developing map due to
R. Bryant \cite[pp. 333-4]{Br87},  Umehara-Yamada \cite[p.
76]{UY00} and Eremenko \cite[p. 3350]{Er04}. Let $g$ be a
conformal metric of constant curvature one on $\Sigma$
representing the divisor $D$. Denote
$\Sigma^*=\Sigma\backslash\{p_1,\cdots,p_n\}$. Every point $p$ in
$\Sigma^*$ has a neighborhood $U_p$ isometric (so conformal) to an
open set ${\frak U}_p$ of the Riemann sphere ${\overline {\Bbb
C}}$ endowed with the standard metric $g_{\rm st}$. Denoting  by
${\frak f}_p:U_p\to {\frak U}_p$ this isometry (conformal map),
Umehara-Yamada and Eremenko loc cit claimed that ${\frak f}_p$ can
be extended to the whole of $\Sigma^*$ by analytic continuation
such that the extension gives a multi-valued locally univalent
meromorphic function $f$ on $\Sigma^*$, whose monodromy belongs to
the group $\SU$ of orientation-preserving isometries of
${\overline {\Bbb C}}$ (cf Lemma \ref{lem:Er}).
 Hence, the metric $g$ can be thought of as the pullback
\[g=\frac{4|f'(z)|^2\,|dz|^2}{\bigl(1+|f(z)|^2\bigr)^2}\]
under $f$ of $g_{\rm st}$. Moreover, stimulated by Umehara-Yamada
\cite[(2.10)]{UY00} and Eremenko \cite[(2)]{Er04}, we show in
Lemma \ref{lem:Sch} that the Schwarzian $\{f,\,z\}$ of $f$ in a
neighborhood $U_j$ of $p_j$ with the complex coordinate $z$ with
$z(p_j)=0$ has the expression
\[\{f,\,z\}:=\frac{f'''(z)}{f'(z)}-
\frac{3}{2}\,\left(\frac{f''(z)}{f'(z)}\right)^2=
\frac{c_j}{z^2}+\frac{d_j}{z}+\psi_j(z),\] where
$c_j=(1-\alpha_j^2)/2$, $d_j$ are constants and $\psi_j$ are
holomorphic functions in $U_j$, dependent on the complex
coordinate $z$. Since the value of $c_j$ is independent of the
choice of the complex coordinate $z$, we say that $f$ is {\it
compatible with the divisor $D=\sum_j (\alpha_j-1)\,P_j$}. We
now arrive at\\

\nd {\bf Definition 1.1} (Umehara-Yamada \cite[p. 76]{UY00}) Let
$g$ be a conformal metric on $\Sigma$ of constant curvature one
and representing the divisor $D$. We call a multi-valued locally
univalent meromorphic function $f$ on $\Sigma^*$ a {\it developing
map} of the metric $g$ if
$g=f^*\,g_{\rm st}$.\\

On the other hand, if there exists a multi-valued meromorphic
function $f$ on $\Sigma^*$, which is compatible with the divisor
$D$ and has monodromy in $\SU$, then there exists a conformal
metric $g=f^*g_{\rm st}$ with constant curvature one and
representing $D$ (cf Lemma \ref{lem:char}). Therefore, we can sum
up the above into a necessary and sufficient condition (cf Theorem
\ref{thm:ns}) for the existence problem of conical conformal
metrics of constant curvature one on $\Sigma$. \\

% Add relations between abelian metrics and HCMU metrics.

In this manuscript, we mainly focus on a special class of conical
conformal metrics of constant curvature one, called {\it reducible
metrics}, which we can classified by using abelian differentials
of 3rd kind.\\

\nd {\bf Definition 1.2} (Umehara-Yamada \cite[p. 76]{UY00})  We
call a conformal metric $g$ on $\Sigma$ of constant curvature one
and with finite conical singularties an {\it irreducible metric}
if the monodromy group of a developing map of the metric $g$ can
not be diagonalized, i.e. the monodromy group has no fixed point
on the Riemann sphere $\BCB$ (cf Lemma \ref{lem:diag}). We call
$g$ {\it reducible} if the monodromy group has at least one fixed
point on $\BCB$.  We call a reducible metric {\it (non-)trivial},
if the monodromy of a developing map of the metric is
(non-)trivial. Lemma \ref{lem:moduli} tells us that these
definitions do not depend on the choice of a developing map.\\

A trivial reducible metric is a pullback of $g_{\rm st}$ under
some rational function on $\Sigma$ (Lemma \ref{lem:trivial}). Each
subgroup of $\SU$ having at least one fixed point on $\BCB$ is
abelian, and, up to conjugacy, it can be thought of as a subgroup
of the standard maximal torus
$${\rm
U}(1)=\left\{{\rm diag}\,\Bigl(e^{\sqrt{-1}\theta},\,
e^{-\sqrt{-1}\theta}\Bigr):\,\theta\in {\Bbb R}\right\}/\{\pm
I_2\}$$ of $\SU$ (cf Lemma \ref{lem:diag}). Each fractional
transformation in $\U$ is the multiplication by
$e^{2\sqrt{-1}\theta}$. Therefore, for a non-trivial reducible
metric $g$ on $\Sigma$,  by Lemma \ref{lem:moduli} there exists
exactly two developing maps, say $f$ and $1/f$, of the metric
$g$, whose monodromies belong to $\U$.\\

\nd {\bf Definition 1.3} Let $g$ be a non-trivial reducible metric
on a compact Riemann surface $\Sigma$. We call a developing map
$f$ of $g$ {\it multiplicative} if the monodromy of $f$ belongs to
$\U$. Such $f$ is unique up to reciprocal, and the logarithmic
differential
$$\omega:=d\,\bigl(\log\,f\bigr)=\frac{df}{f}$$ of the
multiplicative developing map $f$ is a meromorphic 1-form on
$\Sigma^*$. Actually, $\omega$ can be extended to be an abelian
differential of 3rd kind on $\Sigma$ (cf Lemma \ref{lem:abelPsi}),
which we call the {\it character 1-form} of the reducible metric
$g$. Hence the character 1-form of a non-trivial reducible metric
is unique up to sign.

Let $g$ be a trivial reducible metric on $\Sigma$. By Lemma
\ref{lem:trivial}, there exists a rational function $f:\Sigma\to
\BCB$ such that $g=f^*g_{\rm st}$. By Lemma 1.1,  each developing
map of the metric $g$ is a rational function and multiplicative.
We call by a {\it character 1-form} of the metric $g$ the
logarithmic differential of a developing map of $g$. The character
1-forms of the trivial reducible metric $g$ are automatically
abelian differentials of 3rd kind on $\Sigma$. \\

To set the notations for stating the properties of character
1-forms, we need say something more about the standard metric
$g_{\rm st}$ on the Riemann sphere $\BCB$, which is a trivial
reducible metric. The set of all developing maps of $g_{\rm st}$
can be identified with the group $\SU$. Up to reciprocal, any two
developing maps of $g_{\rm st}$, fixing $0$ and $\infty$,
respectively, differ by a multiple complex constant with modulus
1. Up to sign, the logarithmic differentials of all the developing
maps, fixing $0$ and $\infty$, respectively, coincide with the
abelian differential
$\Theta:=d\,\bigl(\log\,w\bigr)=\frac{dw}{w}$, which has two
simple poles of $0$ and $\infty$. The residues of $\Theta$ at $0$
and $\infty$ equal $1$ and $-1$, respectively. The algebraic dual
 $X:=w\,\frac{\partial}{\partial w}$ of $\Theta$ is a meromorphic
vector field with two simple zeroes of $0$ and $\infty$. The index
of $X$ equals $1$ at both $0$ and $\infty$.
$\Phi(w)=\frac{4|w|^2}{1+|w|^2}$ is a smooth Morse function on
$\BCB$, whose complex gradient field $\Phi^{,\,w}\,
\frac{\partial}{\partial w}$ equals $X$. Moreover, $\Phi$ has only
two critical points, which are the minimal point $0$ and the
maximal point $\infty$. Consider a multiplicative developing map
$f$ of a reducible metric $g$ on $\Sigma$.  Up to sign, the
character 1-form $\omega=\frac{df}{f}=d\,\bigl(\log\,f\bigr)$
equals the pullback $f^*\,\Theta$ of $\Theta$ by $f$. Denote by
$Y:=\frac{f(z)}{f'(z)}\, \frac{\partial}{\partial z}$ the
algebraic dual vector field of $\omega$, which is a meromorphic
vector field on $\Sigma$. Then $Y$ equals the complex gradient
field $\Psi^{,\,z}\, \frac{\partial}{\partial z}$ of the smooth
function $\Psi(z)=\frac{4|f(z)|^2}{1+|f(z)|^2}$ on $\Sigma^*$,
which can be continuously extended to $\Sigma$ (cf Lemma
\ref{lem:abelPsi}).

Using the above notations, we state more precisely
the properties of the character 1-form of a reducible metric. \\

%%%%%%%%%%%%% dual vector field of the character 1-form

\begin{thm}
\label{thm:nontrivial} Let $g$ be a reducible metric representing
the divisor $D=\sum_{j=1}^n\,(\alpha_j-1)\,P_j$ with
$1\not=\alpha_j>0$. Let $f$ be a developing map of $g$ and be
multiplicative if $g$ is non-trivial such that the character
1-form $\omega=\frac{df}{f}$ of $g$ equals $f^*\Theta$. Let $Y$ be
the algebraic dual vector field of $\omega$, and
$\Psi(z)=\frac{4|f(z)|^2}{1+|f(z)|^2}$. Then the
following statements hold.\\

\nd {\bf (1)} The set of zeroes of the meromorphic vector field
$Y$ coincides with the extremal point set of the function $\Psi$.
Each zero of $Y$ is simple, and $Y$ vanishes at each point $p_j$
where $\alpha_j>0$ is a non-integer. The set of poles of $Y$
coincides with the saddle point set of $\Psi$. Each pole of $Y$ is
some conical singularity $p_j$ of the reducible metric $g$, where
$\alpha_j$ is an integer greater than 1 and the order of the pole
$p_j$ of $Y$ equals $\alpha_j-1$.\\

\nd {\bf (2)}  Let $p_1,\cdots,p_J$ be the saddle points of
$\Psi$, $p_{J+1},\cdots,p_n$ the singular extremal points of
$\Psi$, and $e_1,\cdots,e_S$ the smooth extremal points of $\Psi$
on $\Sigma^*$. Then the canonical divisor of the character 1-form
$\omega$ has form
\[(\omega)=\sum_{j=1}^J\,(\alpha_j-1)\,P_j-\sum_{k=J+1}^n\,P_k-\sum_{\ell=1}^S\,
E_\ell.\] In particular, each pole of $\omega$ is simple, i.e.
$\omega$ is an abelian differential of 3rd kind. The residue of
$\omega$ at the pole $e_\ell$ equals $1$ or $-1$, where $e_\ell$
is a minimal or maximal point of $\Psi$; the residue of $\omega$
at the pole $p_k$ equals $\alpha_j$ or $-\alpha_j$, where $p_k$ is
a minimal or maximal point of $\Psi$. Moreover, the real part of
$\omega$ is exact on $\Sigma':=\Sigma\backslash\{p_{J+1},\cdots,
p_n,e_1,\cdots, e_S\}\supset\Sigma^*$,
\[2\,\Re\,\omega=d\,(\log\,|f|^2).\]

\nd {\bf (3)} The multiplicative developing map $f$ on
$\Sigma'\supset \Sigma^*$ can be expressed by
\[f(z)=\exp\,\Bigl(\int^z\, \omega\Bigr).\]
In particular, the local monodromy of $f$ around each $p_j$ {\rm (
$1\leq j\leq J$)} is trivial, and the limit $\lim_{p\to p_j}\,
f(p)$ exists and belongs to ${\Bbb C}\backslash\{0\}$.  If we
continue analytically a function element ${\frak f}$ along a
simple and sufficiently small loop winding around $p_k$ {\rm
($J+1\leq k\leq n$)} on counterclockwise, then we obtain  ${\frak
f}\,\exp\,\bigl(2\pi\sqrt{-1}\alpha_k\bigr)$. The limit
$\lim_{p\to p_k}\, |f(p)|$ exists, and equals $0$ or $+\infty$,
provided $p_k$ is a minimal or maximal point of $\Psi(z)$.
It is also the case for $e_\ell$. \\

%\nd {\bf (4)} There are a finite number of geodesics which connect
%extremal points and saddle points of $\Psi$, and $\Sigma$ can be
%decomposed into a finite number pieces by cutting along these
%geodesics where each piece is a geodesic bigon in the Riemann
%sphere $\BCB$.

\end{thm}

Using abelian differentials of 3rd kind with the above properties,
we can construct new examples of conformal metrics with constant
curvature one and with finite conical singularities. \\

\begin{thm}
\label{thm:char} Let $\omega$ be an abelian differential of 3rd
kind having poles on a compact Riemann surface $\Sigma$, whose
residues are all nonzero real numbers and whose real part is exact
outside the set of poles of $\omega$. Then there exists a unique
reducible metric $g$ on $\Sigma$ such that $\omega$ is one of the
character 1-forms of $g$ and $g$ can be expressed by $g=f^*\,
g_{\rm st}$, where
\[f(z)=\exp\,\Bigl(\int^z\, \omega\Bigr)\]
is a multi-valued locally univalent meromorphic function on
$\Sigma\backslash\{\rm{poles\ of}\ \omega\}$ with monodromy in
$\U$. Suppose that the canonical divisor of $\omega$ has form
\[(\omega)=\sum_{j=1}^J\,(\alpha_j-1)\,P_j-\sum_{k=J+1}^N\,Q_k,\]
where $\alpha_j$ are integers $>1$. Then the divisor $D$
represented by $g$ has form
\[D=\sum_{j=1}^J\,(\alpha_j-1)\,P_j+
\sum_{k=J+1}^N\,\Bigl(|{\rm Res}_{Q_k}(\omega)|-1\Bigr)\, Q_k\]

Moreover, $g$ is a trivial reducible metric if and only if the
integral of $\omega$ on each loop in $\Sigma\backslash\{\rm{poles\
of}\ \omega\}$ are $2\pi\sqrt{-1}$ times integers. In particular,
each residue of $\omega$ is an integer.
\end{thm}

\nd{\bf Remark 1.1.} Each reducible metric does not satisfy
Troyanov's condition \eqref{equ:tro} (cf Corollary \ref{cor:tro}).
Therefore, we obtain a class of new examples of conformal metrics
of constant curvature one with finite singularities, since there
exist plenty of abelian differentials of 3rd kind satisfying the
condition in Theorem \ref{thm:char} (cf Springer \cite[Corollary
8-3]{Spr57}).

Troyanov's condition \eqref{equ:tro} for the corresponding
irreducible metrics only depends on the values of angles. Example
\ref{exam:abel} shows that the existence of reducible metrics does
not only depends on angles, but also on the position of
singularities.\\

\nd{\bf Remark 1.2.} Umehara-Yamada \cite{UY00} called a
non-trivial reducible metric ${\cal H}^1$-reducible, and a trivial
reducible metric ${\cal
H}^3$-reducible.\\

\nd {\bf Remark 1.3.} Besides Umehara-Yamada \cite{UY00}, our
motivation of defining character 1-form comes from \cite{CW11,
CWX12}, where the authors use character 1-form to make the
complete classification of HCMU metrics of non-constant curvature
on compact Riemann
surfaces. We say more words about this in the ending of Section \ref{sec:pfthm}.\\

\begin{thm}
\label{thm:rational} A conformal metric of constant curvature one,
representing an effective ${\Bbb Z}$-divisor on the 2-sphere is a
trivial reducible metric. That is, it is the pullback under some
rational function on the 2-sphere of the standard metric $g_{\rm
st}$ on the Riemann sphere $\BCB$.\end{thm}

\nd {\bf Remark 1.4.} The case where $D$ supports at two or three
points was proved in \cite{Tr89, FH, UY00, Er04}. However, Theorem
\ref{thm:rational} does not hold in general for compact Riemann
surfaces of nonzero genus (cf Example \ref{exam:nonabel}).
\\

\nd{\bf Remark 1.5.} Since we only consider conformal metrics with
finite area in this manuscript, the singularities of "zero angle",
would not show up (cf \cite[Proposition 4]{Br87}). We prove a much
more general result in Section 5
than Proposition 4 in \cite{Br87} using a different method. \\

We explain the organization of this paper. In Section 2, we shall
firstly make a detailed exposition on developing map in Lemmata
\ref{lem:Er} and \ref{lem:moduli}. We compute, in Lemma
\ref{lem:Sch} of Section 3, the Schwarzian of a developing map of
a conformal metric with constant curvature one and representing a
divisor $D$.  In Lemma \ref{lem:char} we show that the converse of
Lemma \ref{lem:Sch} also  holds true. Then Theorem
\ref{thm:rational} follows from these two lemata.  In Section 4 we
prove Theorems \ref{thm:nontrivial} and \ref{thm:char} as an
application of Lemma \ref{lem:char}. These two theorems is applied
to giving some examples of irreducible and reducible metrics in
Corollary \ref{cor:tro} and Examples \ref{exam:nonabel} and
\ref{exam:abel}. In Section 5, we give an alternative proof to a
theorem of Troyanov \cite{Tr89}, as an application of Theorems
\ref{thm:rational} and \ref{thm:nontrivial}. Moreover, we discuss
the non-uniqueness of reducible metrics representing a given
divisor $D$. At last, we propose some speculations about both
irreducible and reducible metrics conformal metrics. In the last
section, it is proved that no singularity with angle "zero" would
appear in the local geometry of nonnegative curvature.\\

%Moreover, we shall describe briefly character 1-forms associated
%with the Kahler-Einstein metric on a compact complex manifold with
%constant holomorphic sectional curvature one with conical
%singularities along a normal crossing divisor.

\nd {\bf Acknowledgments.}   The last author would like to express
his deep gratitude to Professor Xiuxiong Chen for his constant
moral support and lots of invaluable conversations. He thanks
Professor Masaaki Umehara so much for many valuable comments and
discussions through e-mails. He is also very grateful to Professor
Ryoichi Kobayashi, Dr Zhi Chen and Dr Jinxing Xu for their
stimulating conversations. In particular, the penetrating
questions from Dr Jinxing Xu give him the impetus to fixing a gap in the old
version of the manuscript.

\section{Existence of developing maps and their monodromy}

\begin{lem}
\label{lem:Er} Let $g$ be a conformal metric on a compact Riemann
surface $\Sigma$ of constant curvature one and representing the
divisor $D=\Sigma_{j=1}^n\, (\alpha_j-1)\,P_j$ with $\alpha_j>0$.
Then there exists a multi-valued locally univalent holomorophic
map $f$ from $\Sigma^*:=\Sigma\backslash\{p_1,\cdots,p_n\}$ to the
Riemann sphere $\BCB$ such that the monodromy of $f$ belongs to
$\SU$ and $$g=f^*\,g_{\rm st}.$$ Recall that $g_{\rm
st}=\df{4|dw|^2}{(1+|w|^2)^2}$ is the standard metric over $\BCB$.
\end{lem}

\nd \pf Denote by $d(\cdot,\,\cdot)$ the distance on $\Sigma$
induced by the metric $g$. Choose an arbitrary point $p$ in
$\Sigma^*$ and fix it. Take a positive number $r=r_p$ sufficiently
small such that $d(p,\,\{p_1,\cdots,p_n\})>r$ and there exists a
geodesic polar coordinate chart in the open metric ball
$B(p,\,r)=\{q\in \Sigma:\,d(p,\,q)<r\}\subset\Sigma^*$. Choose a
positively oriented orthonormal basis
$\{e^1,\,e^2\}=\{e^1_p,\,e^2_p\}$ of the tangent space
$T_p\,\Sigma$. Choose an arbitrary point ${\frak p}\in\BCB$ and
fix it. Since the Gauss curvature of $\left(B(p,\,r),\,g\right)$
is constant one, by a theorem of Riemann (\cite[p.136]{Pe06}),
there exists an open metric ball ${\frak B}({\frak p},\,r)$ in the
Riemann sphere $(\BCB,\,g_{\rm st})$ and an orientation-preserving
isometry ${\frak f}_p$ from $\left(B(p,\,r),\,g\right)$ onto
$\left({\frak B}({\frak p},\,r),\,g_{\rm st}\right)$. Denote
${\frak e}^1_{\frak p}:={\frak f}_*(e^1),\,{\frak e}^2_{\frak
p}:={\frak f}_*(e^2)$, which is also a positively oriented
orthonormal basis of $T_{\frak p}\,\BCB$. Then ${\frak f}_p$ is a
conformal map from $B(p,\,r)$ to ${\frak B}({\frak p},\,r)$.

Take an arbitrary point $q$ in $\Sigma^*$ and a curve
$L:[0,\,1]\to\Sigma^*$ joining $p$ to $q$. Then there exists
$$0<\delta<\min\left(d\left(\gamma([0,\,1]),\,
\{p_1,\cdots,p_n\}\right),\,r_p\right)$$ such that  there exists a
geodesic polar coordinate chart in the open metric ball
$B(a,\,\delta)$ for each point $a$ on the curve $L$. If we divide
properly the interval $0\leq t\leq 1$ into $n$ subintervals for
sufficiently large $n$, $0=\gamma_0<\gamma_1<\cdots<\gamma_n=1$,
then the curve $L$ splits into $n$ subarcs $L_1,\,L_2,\,\cdots,
L_n$ with $L_g$ ($g=1,\cdots,n$) joining
$L(\gamma_{g-1})=:c_{g-1}$ to $L(\gamma_g)=:c_g$. Moreover, if we
denote by $B_0,\,B_1,\,\cdots,\,B_n$ the open metric balls with
centers $a=c_0,\,c_1,\,\cdots,\,c_n$ and with radius $\delta$,
then the closed arcs $L_g$ lie completely in $B_{g-1}$ for
$g=1,\cdots,n$. Let $f_0$ be the restriction to $B_0$ of the
conformal map ${\frak f}_p:B(p,\,r_p)\to {\frak B}({\frak
p},\,r_p)$. Then $f_0$ is an isometry (conformal map) from $B_0$
onto ${\frak B}_0:={\frak B}({\frak p},\,\delta)$. Choose a
positively oriented orthonormal basis $\{e^1_{c_1},\,e^2_{c_2}\}$
of $T_{c_1}\,\Sigma$. Since $c_1\in B_0$, we denote ${\frak
c}_1:=f_0(c_1)\in {\frak B}_0$, ${\frak
e}^1_{c_1}:=(f_0)_*(e^1_{c_1})$ and ${\frak
e}^2_{c_1}:=(f_0)_*(e^2_{c_1})$. Then there exists a unique
isometry $f_1:B_1\to {\frak B}_1:={\frak B}({\frak c}_1,\,\delta)$
such that $f_1(c_1)={\frak c}_1$ and $f_1$ maps
$\{e^1_{c_1},\,e^2_{c_1}\}$ to $\{{\frak e}^1_{{\frak
c}_1},\,{\frak e}^2_{{\frak c}_1}\}$. Then $f_1=f_0$ on $B_0\cap
B_1$. Since $c_2\in L_2\subset B_1$, $f_1$ is an analytic
continuation of $f_0$ from point $c_0$ to $c_2$ along the arc
$L_0\cup L_1$. In this way, we obtain $f_0,\,\cdots,\,f_n$, which
are recursively defined on $B_0,\,\cdots,B_n$ and give an analytic
continuation of ${\frak f}_p$ from $p$ to $q$ along the curve $L$.
Using the same argument as \cite[pp.13-15]{Si69}, we can show that
this analytic continuation is independent of the choice of
division points on $L$. Moreover, if $L^*$ is another curve in
$\Sigma^*$ joining $a$ to $b$ and homotopic to $L$, then the
result of doing analytic continuation of ${\frak f}_p$ along $L^*$
is the same as along $L$. Summing up, we obtain a multi-valued
locally isometrical (univalent conformal) map $f$ from
$(\Sigma^*,\,g)$ to $(\BCB,\,g_{\rm st})$.

At last, we prove that all the monodromy of $f$ belongs to $\SU$.
Let $L:[0,\,1]\to \Sigma^*$ be a closed curve with $L(0)=L(1)=p$.
We use the notation in the previous paragraph. Recall that $f_0$
maps $p=c_0$ to ${\frak p}$ and $(f_0)_*$ maps $\{e^1,\,e^2\}$ to
$\{{\frak e}^1_{\frak p},\,{\frak e}^2_{\frak p}\}$, and $f_n$
maps $p=c_n$ to ${\frak c}_n$ and $(f_n)_*$ maps $\{e^1,\,e^2\}$
to $\{{\frak e}^1_{{\frak c}_n},\,{\frak e}^2_{{\frak c}_n}\}$.
Then there exists a unique isometry ${\frak L}\in \SU$ of
$(\BCB,\,g_{\rm st})$ such that ${\frak L}({\frak p})={\frak c}_n$
and ${\frak L}_*$ maps $\{{\frak e}^1_{\frak p},\,{\frak
e}^2_{\frak p}\}$ to $\{{\frak e}^1_{{\frak c}_n},\,{\frak
e}^2_{{\frak c}_n}\}$. Therefore $f_n={\frak L}\circ f_0$. \hf\\

\nd \begin{lem} \label{lem:moduli}  Any two developing maps
$f_1,\,f_2$ of the metric $g$ are related by a fractional linear
transformation ${\frak L}\in\SU$, $f_2={\frak L}\circ f_1.$ In
particular, any two developing maps of $g$ have mutually conjugate
monodromy in $\SU$. Then we call this conjugate class the
monodromy of the metric $g$. The space of developing maps of the
metric $g$ has a one-to-one correspondence with the quotient group
of $\SU$ by the monodromy group of a developing map of $g$.
\end{lem}

\nd\pf  Take a point $p\in\Sigma^*$ and a positively oriented
orthonormal basis $\{e^1,\,e^2\}$ of $T_p\,\Sigma^*$. Let ${\frak
f}_j$ be a function element of $f_j$ near $p$ for $j=1,\,2$.
Denote ${\frak p}_j:={\frak f}_j(p)$ and ${\frak e}^k_{{\frak
p}_j}:=({\frak f}_j)_*(e^k)$ for $j,\,k=1,\,2$. Then there exists
a unique ${\frak L}\in\SU$ such that ${\frak L}({\frak
p}_1)={\frak p}_2$ and ${\frak L}_*$ maps $\{{\frak e}^1_{{\frak
p}_1},\,{\frak e}^2_{{\frak p}_1}\}$ to $\{{\frak e}^1_{{\frak
p}_2},\,{\frak e}^2_{{\frak p}_2}\}$. Then we obtain the equality
${\frak f}_2={\frak L}\circ {\frak f}_1$ near $p$, which implies
$f_2={\frak L}\circ f_1$. It follows from direct computation that
the monodromy of $f_1$ and $f_2$ are mutually conjugate.

Given a developing map $f$ and a fractional linear transformation
${\frak L}\in \SU$, we can see ${\frak L}\circ f=f$ if and only if
there exists a point $p\in\Sigma^*$ and a functional element
${\frak f}_p$ near $p$ such that ${\frak L}\circ {\frak f}_p$ is
another function element of $f$ near $p$. That is, ${\frak L}$
belongs to the image of the monodromy representation of
$\pi_1(\Sigma^*,\, p)$ with respect to ${\frak f}$.
 Therefore, $\SU$ acts in
this way transitively on the set of all developing maps with
isotropy group isomorphic to the monodromy group.
\hf\\

\nd {\bf Remark 2.1.}  The developing maps also exist for flat or
hyperbolic conformal metrics with finite conical or cusp
singularities, and an analogue of Lemmata \ref{lem:moduli} and
\ref{lem:char} holds.
\\

\section{The Schwarzian of a developing map}
\begin{lem}
\label{lem:Sch} Let $g$ be a conformal metric of constant
curvature one on a compact Riemann surface $\Sigma$ and $g$
represent a divisor $D=\Sigma_{j=1}^n\, (\alpha_j-1)\,P_j$, where
$\alpha_j>0$ for all $j$. Suppose that
$f:\Sigma^*=\Sigma\backslash \{p_1,\cdots,p_n\}\to \BCB$ is a
developing map of $g$. Then the Schwarzian $\{f,\,z\}$ of $f$
equals
\[\{f,\,z\}=
\frac{1-\alpha_j^2}{2\,z^2}+\frac{d_j}{z}+\psi_j(z)\] in a
neighborhood $U_j$ of $p_j$ with the complex coordinate $z$ and
$z(p_j)=0$, where $d_j$ are constants and $\psi_j$ are holomorphic
functions in $U_j$, dependending on the complex coordinate $z$.
\end{lem}

\nd \pf If we rewrite the metric
$g=\frac{4|f'(z)|^2\,|dz|^2}{\bigl(1+|f(z)|^2\bigr)^2}$ as
$g=e^{2u}\,|dz|^2$, then we find $u=\log\,|f'(z)|+\log\,
2-\log\,(1+|f|^2)$. The lemma in \cite[p.300]{Tr89} tells us that
\[\eta(z)=2\left(\frac{\partial^2 u}{\partial
z^2}-\left(\frac{\partial u}{\partial z}\right)^2\right)\, dz^2\]
defines a {\it projective connection} compatible with the divisor
$D$. The interested reader could find in \cite{Tr89} the
definition of the projective connection, which we will not use in
this paper, however.  The compatibility of the projective
connection $\eta$ with the divisor $D$ in \cite[p.300]{Tr89} means
that
\[\eta(z)=\left(\frac{1-\alpha_j^2}{2\,z^2}+\frac{d_j}{z}+\phi_j(z)\right)\,dz^2,
\quad \phi_j\quad \text{holomorphic},\] where $z$ is the complex
coordinate near $p_j$. Since the developing map $f$ is a
projective multi-valued function on $\Sigma^*$, its Schwarzian
$\{f,\,z\}$ with respect to the complex coordinate $z$ near $p_j$
is a single valued function of $z$. At last, we find
\begin{eqnarray*}
2\left(\frac{\partial^2 u}{\partial z^2}-\left(\frac{\partial
u}{\partial z}\right)^2\right) &=&
2\tpds{}{z}\left(\df{f"(z)}{2f'(z)}-\df{f'(z)\overline{f}}{1+|f|^2}\right)
-2\left(\df{f"(z)}{2f'(z)}-\df{f'(z)\overline{f}}{1+|f|^2}\right)^2\\
&=&\left(\df{f'''(z)}{f'(z)}-\left(\df{f"(z)}{f'(z)}\right)^2-
\df{2f"(z)\overline{f(z)}}{1+|f|^2}+
2\left(\df{f'(z)\overline{f}}{1+|f|^2}\right)^2\right)\\
&& -\left(\df{1}{2}\left(\df{f"(z)}{f'(z)}\right)^2-
\df{2f"(z)\overline{f(z)}}{1+|f|^2}+
2\left(\df{f'(z)\overline{f}}{1+|f|^2}\right)^2\right)\\
&=&\df{f'''(z)}{f'(z)}-\df{3}{2}\left(\df{f"(z)}{f'(z)}\right)^2=\{f,z\}.
\end{eqnarray*}
\hf\\

A multi-valued locally univalent meromorphic function $h$ on
$\Sigma^*$ is said to be {\it projective} if any two function
elements ${\frak h}_1,\, {\frak h_2}$ of $h$ near a point
$p\in\Sigma^*$ are related by a fractional linear transformation
$L\in \PGL$,
${\frak h}_1=L\circ{\frak h}_2$.  \\

%%%%%%%%%%%%% Construct a csc metric with a multi-valued functon

\nd \begin{lem} \label{lem:char} Let $f:\Sigma^*\to \overline{\Bbb
C}$ be a projective multi-valued locally univalent meromorphic
function, and the monodromy of $f$ belongs to a maximal compact
subgroup of ${\rm PGL}(2,\,{\Bbb C})$. If $f$ has regular
singularity of weight $0\not=c_j<1/2$ at $p_j$ for all
$j=1,\cdots, n$, then there exists a neighborhood $U_j$ of $p_j$
with complex coordinate $z$ and $L_j\in \PGL$ such that $z(p_j)=0$
and $g_j=L_j\circ f$ has the form $g_j(z)=z^{\alpha_j}$, where
$0<\alpha_j\not=1$ and $c_j=(1-\alpha_j^2)/2$. Moreover, there
exists ${\frak L}\in \PGL$ such that the pullback $({\frak L}\circ
f)^*\,g_{\rm st}$ of the standard metric $g_{\rm st}$ by ${\frak
L}\circ f$ is a conformal metric of constant curvature one, which
represents the divisor $D=\sum_j\,(\alpha_j-1)\, P_j$. In
particular, if the monodromy of $f$ belongs to $\SU$, then the
fractional linear transformation ${\frak L}$ turns out to be the
identity map.
\end{lem}

\nd \pf \quad Recall the well known fact that every maximal
compact group of $\PGL$ is conjugate to the subgroup $\SU$. There
exists a fractional linear transformation ${\frak L}$ such that
the monodromy of ${\frak L}\circ f$ belongs to $\SU$. Hence, we
may assume that it is the case for $f$ without loss of generality.

We firstly show the first statement of Lemma 1.2 that {\it there
exists a neighborhood $U_j$ of $p_j$ with complex coordinate $z$
and some $L_j\in\PGL$ such that $g_j=L_j\circ f$ has form
$z^{\alpha_j}$}.
Since $f$ has a regular singularity of weight
$c_j<1/2$, we could choose a neighborhood $U_j$ of $p_j$ and a
complex coordinate $x$ on $U_j$ such that $x(p_j)=0$ and
$$\{f,\,x\}=\frac{c_j}{x^2}+\frac{d_j}{x}+\phi_j(x),$$ where
$\phi_j(x)$ is holomorphic in $U_j$. By \cite[p.39,
Proposition]{Yo87}, in the neighborhood $U_j$, there are two
linearly independent solutions $u_0$ and $u_1$ of the equation
\[\frac{d^2 u}{dx^2}+
\frac{1}{2}\left(\frac{c_j}{x^2}+\frac{d_j}{x}+\phi_j(x)\right)u=0\]
with single valued coefficient such that
$f(x)=\frac{u_1(x)}{u_0(x)}$. Actually, we have
$u_0=\left(\frac{df}{dx}\right)^{-1/2}$ and $u_1=f(x)\,u_0$.
Moreover, if $f$ changes projectively \[f\mapsto
\frac{af+b}{cf+d}\quad \text{with}\quad ad-bc=1,\] then $u_0$ and
$u_1$ change linearly
\[\begin{pmatrix}
u_0\\
u_1
\end{pmatrix}
 \mapsto \begin{pmatrix} d & c\\
b & a
\end{pmatrix}
\begin{pmatrix}
u_0\\
u_1
\end{pmatrix}, \]
and vice versa.

Consider the operator $L_j:=x^2\frac{d^2}{dx^2}+q_j(x)$ with
$q_j(x)=\left(c_j+d_j x+x^2\,\phi_j(x)\right)/2$. Then both $u_0$
and $u_1$ are solutions of the equation $L_j\, u=0$. Since the
monodromy of $f$ belongs to $\SU$, the cyclic group generated by
the local monodromy of the equation $L_j\,u=0$ around the $x=0$ is
contained in a maximal compact subgroup of $\PGL$ conjugate to
$\SU$. Note that the equation $L_j\,u=0$ has regular singularity
at $0$. We could apply the Frobenius method (cf \cite[$\S\,$
2.5]{Yo87}) to solving it. Choose $\alpha_j>0$ such that
$c_j=(1-\alpha_j^2)/2$. Then the indicial equation
\begin{equation}
\label{ind} f(s)=s(s-1)+\frac{1-\alpha_j^2}{4}=0
\end{equation} of the
differential equation $L_j u=0$ at $x=0$ has roots
$s_0=\frac{1-\alpha_j}{2}$ and $s_1=\frac{1+\alpha_j}{2}$, and
$s_1-s_0=\alpha_j>0$. Let $\sum_{k=0}^\infty\, b_k\,x^k$ be the
power series expansion of $q_j(x)$, where $b_0=c_j/2$. Let $s$ be
a parameter. Then $u(s,x)=x^s\,\sum_{k=0}^\infty\,c_k(s)\,x^k$
with $c_0(s)\equiv1$ is an solution of $L_ju=0$ if and only if the
equation
\begin{equation}
f(s+n)\,c_n+R_n=0 \tag{$\sharp_n$}
\end{equation}
holds for all $n=0,1,2,\cdots$, where
\begin{eqnarray*}
R_0=0,\quad \text{and for}\quad n>0,\quad
R_n=R_n(c_1,\cdots,c_{n-1},s)=\sum_{i=0}^{n-1}\, c_i\,b_{n-i}.
\end{eqnarray*}
Note that the equation ($\sharp_0$) is exactly the indicial
equation \eqref{ind}.  Since $f(s_1+n)\not=0$ for all $n\geq 1$,
we find that $u(s_1,\,x)$ is a
solution of the equation. \\

{\it Case 1}\quad Suppose that $s_1-s_0=\alpha_j$ is not an
integer. Then by the same reason $u(s_0, \,x)$ is another
solution, which is linearly independent of $u(s_1,\,x)$. Summing
up, we have
\[u(s_0,\,x)=x^{s_0}\,\left(1+\psi_0\right)\quad \text{and}\quad
 u(s_1,\,x)=x^{s_1}\,\left(1+\psi_1\right),\]
where both $\psi_0$ and $\psi_1$ are holomorophic functions
vanishing at $0$. Here we take a smaller neighborhood of $0$ than
$U_j$ to assure the convergence of the power series defining
$\psi_k$ if necessary. Since both $u_0(x)$ and $u_1(x)$ are linear
combinations of $u(s_0,\,x)$ and $u(s_1,\,x)$,
$f(x)=u_1(x)/u_0(x)$ equals some fractional linear transform of
$u(s_1,\,x)/u(s_0,\,x)$. For simplicity of notation, we may assume
$f(x)=u(s_1,\,x)/u(s_0,\,x)$ equals $x^{\alpha_j}$ times a
holomorphic function $\varphi_j(x)$ with $\varphi_j(0)=1$.
Therefore, we could choose another complex coordinate $z=z(x)$ of
$U_j$,
under which $f=f(z)=z^{\alpha_j}$. \\

{\it Case 2} \quad Suppose that $m:=s_1-s_0=\alpha_j$ is an
integer $\geq 2$.

{\it Subcase 2.1}\quad If $R_m=0$, we can solve the equation
($\sharp_n$) for $s=s_0$ for all $n\geq 1$ by choosing $c_m$
arbitrarily, and obtain another solution $u(s_0,\,x)$ linearly
independent of $u(s_1,\,x)$. The similar argument as Case 1
completes the proof.

{\it Subcase 2.2}\quad Suppose $R_m\not=0$. Define
\[u^*=x^{s_0}\,\sum_{k=0}^\infty\, c_k(s_0)\,x^k,\]
where $c_0=1$, the $c_j$'s ($1\leq j<m$) are determined by
$(\sharp_j)$, while $c_m$ is arbitrarily fixed, and the $c_j$'s
($j>m$) are determined also by $(\sharp_j)$. Then the linear
combination of $u^*$ and $\tpds {}{s}\,u(s,\,x)|_{s=s_1}$
\[U_0(x):=f'(s_1)\,u^*-R_m\,\tpds {}{s}\,u(s,\,x)|_{s=s_1}\]
is a solution. It should be mentioned that since $f$ and $R_n$ are
holomorphic with respect to $s$, it is also the case for the
$c_n$'s. Then we correct a typo in \cite[p.23]{Yo87} and find the
two linearly independent solutions given by
\[
\begin{pmatrix}
U_0(x)\\
u(s_1,\,x)
\end{pmatrix}
=\begin{pmatrix}
x^{s_0} & x^{s_1}\,\log\,x\\
0 & x^{s_1}
\end{pmatrix}
\cdot
\begin{pmatrix}
f'(s_1)\,\sum_{k=0}^\infty\, c_k(s_0)\, x^k-
R_m\,x^m\,\sum_{k=0}^\infty\, c'_k(s_1)\, x^k\\
\sum_{k=0}^\infty\, c_k(s_1)\,x^k
\end{pmatrix}.
\]
Then the local monodromy of the
equation $L_j\,u=0$ at $x=0$ is the conjugacy class in $\PGL$ of
the matrix
$$M=\begin{pmatrix}
1 & 2\pi\sqrt{-1}\\
0 & 1
\end{pmatrix}.$$
However, the cyclic group generated by $M$ is a free abelian
group, which has no limit point under the usual topology of
$\PGL$. Contradict the fact that the monodromy of the equation
$L_j\,u=0$ is contained in a compact group of $\PGL$. That is, we
rule out Subcase 2.2.

Summing up, we prove the statement where $\alpha_j$ is an integer
$\geq 2$. Moreover, we can also see that in this case the local
monodromy at $p_j$ is trivial, i.e. $p_j$ is an apparent
singularity of the equation $L_j\,u=0$ and the multi-valued
function $f$.\\

Since $f$ is locally univalent on $\Sigma^*$ and has monodromy
belonging to $\SU$, $f^*\,g_{\rm st}$ is a well defined smooth
Riemannian metric on $\Sigma^*$ with constant curvature one. The
first statement proved just now implies that this metric has
conical singularities at $p_j$ with angles $2\pi\,\alpha_j$.\hf\\

\nd{\bf Remark 3.1.} Lemma \ref{lem:char} has some overlapping
with Proposition 4 in Bryant \cite{Br87} in the sense that both of them
say the same thing near each singularity. \\

We sum up the above two lemmata into

\begin{thm}
\label{thm:ns} There exists a conformal metric of constant
curvature one and representing a divisor $D$ on a compact Riemann
surface $\Sigma$ if and only if there is a projective multi-valued
meromorphic function on $\Sigma^*=\Sigma\backslash {\rm Supp}\, D$
compatible with the divisor $D$ and having monodromy in $\SU$.
\end{thm}

\nd{\bf Proof of Theorem 1.1}\quad Let $f$ be a developing map of
the metric $g$ representing an effective ${\Bbb Z}$-divisor
$\sum_j\, n_j\,P_j$ on the sphere. By Lemma \ref{lem:Sch}, $f$ has
regular singularity of weight $(1-n_j^2)/2$ at $p_j$. By Lemma
\ref{lem:Er}, the monodromy of $f$ belongs to $\SU$. By Lemma
\ref{lem:char}, there exists ${\frak L}_j\in \PGL$ and a complex
coordinate $z$ near $p$ such that ${\frak L}_j\circ f$ has the
form $f(z)=z^{n_j+1}$ near $p_j$, which implies the local
monodromy of $f$ at $p_j$ is trivial. Since the sphere is simply
connected, the monodromy of $f$ is trivial, i.e. $f$ is a single
valued meromorphic function outside $\{p_j\}$. Morever, $f$ can be
extended meromorphically onto the whole sphere, i.e. $f$ is a
rational function on the sphere.
\hf\\

At a point $p\in \Sigma$ near which $f=f(z)$ is univalent
holomorphic, we find that the Schwarzian $\{f,\,z\}$ is
holomorphic, and vice versa (cf \cite[Remark, p.44]{Yo87}).
Actually, we can prove a
more general result.\\

\begin{lem}
\label{lem:2pi}
 Let $U$ be an open disk containing $0$ in the
complex plane ${\Bbb C}$ with coordinate $w$ and $f$ a projective
multi-valued meormorphic function on $U\backslash\{0\}$ with
regular singularity of weight zero at $0$. That is, $\{f,\,w\}$
equals $\df{d}{w}$ plus a holomorphic function $\phi(w)$, where
both the constant $d$ and $\phi(w)$ depend on the coordinate $w$.
Assume that the subgroup of $\PGL$ generated by the local
monodromy of $f$ at $0$ is precompact in $\PGL$. Then there exists
${\frak L}\in \PGL$ and another complex coordinate $z$ of $U$ such
that ${\frak L}\circ f(z)=z$ and $z(0)=0$. \end{lem}

\nd \pf\quad Use the same argument as Case 2 of the proof of Lemma
1.2. Also note that the indicial equation here has two roots $0$
and $1$. \hf\\

This lemma has the following geometric consequence \\

\nd \begin{prop} \label{prop:2pi}
 The conic singularity with angle
$2\pi$ of a conformal metric with constant curvature one is
actually a smooth point of the metric. \end{prop}

%\begin{lem} \label{lem:cusp} Let $U$ be an open disk containing $0$ in the complex plane ${\Bbb C}$ with coordinate $w$ and $f$ a
%projective multi-valued meormorphic function on $U\backslash\{0\}$ with regular singularity of weight $\df{1}{2}$ at $0$. That is,
%$\{f,\,w\}$ equals $\df{1}{2\,w^2}+\df{d}{w}$ plus a holomorphic function $\phi(w)$, where both the constant $d$ and $\phi$ depend
%on the coordinate $w$. Then the subgroup of $\PGL$ generated by the local
%monodromy at $0$ is not precompact in $\PGL$.
%\end{lem}

%\nd \pf The indicial eqaution here has double root $1/2$, which
%implies that $f$ must have a logarithmic singularity at $0$ (cf
%\cite[p.22]{Yo87}). \hf\\

%Lemma \ref{lem:cusp} has a geometric application as\\

%\begin{prop}
%There exists no conformal metric of constant curvature one, which
%has a cusp singularity.\end{prop}

%\nd \pf (By contradiction) Suppose that there exists such a metric
%with a cusp point $p$. We only need show that near $p$, there
%exists a complex coordinate $z$ such that $z(p)=0$ and the
%developing map $f$ has Schwarzian $\{f,\,z\}$ equals
%$\df{1}{2\,z^2}+\df{d}{z}$ plus a holomorphic function
%$\phi(z)$.\\

%\nd {\bf Remark 4.1.} Wei Wang and Bin Xu prove a more general
%result than Proposition 2.2 in the appendix.

\section{Proof of Theorems 1.2 and 1.3}
\label{sec:pfthm}

\begin{lem}
\label{lem:diag} The following statements hold true.

\nd {\rm (1)} A subgroup $G$ of $\SU$ can be diagonalized if and
only if $G$ has a fixed point on $\BCB$. Such a group is contained
in some maximal torus ${\Bbb T}$ of $\SU$. In particular, $G$ is
abelian.

\nd {\rm (2)} There exists an abelian subgroup of $\SU$ which has
no fixed point on $\BCB$.

\end{lem}

\nd\pf (1) Consider the natural unitary representation $\ro$ of
${\rm SU}(2)$ on $V\cong {\Bbb C}^2$ endowed with the natural
Hermitian inner product $\langle\,,\,\rangle$. It induces a
unitary representation $\ro_H$ on $V$ for each subgroup $H$ of
${\rm SU}(2)$ such that $\ro_H$ is a faithful representation. Let
${\widetilde G}\subset {\rm SU}(2)$ be the lifting of $G\subset
\SU$. We say that $G$ can be {\it diagonalized} if the
representation $(\ro_{\widetilde G},\, V)$ can be decomposed into
the direct sum of two one-dimensional subspaces,
\[V={\Bbb C}e_1\oplus {\Bbb C}e_2\quad{\rm and}\quad
 \langle e_k\,,\,e_\ell\rangle=\delta_{k\ell}.\]
Then, up to conjugacy, ${\widetilde G}$ can be thought of as a
subgroup of the standard maximal torus
$${\rm
U}(1)=\left\{{\rm diag}\,\Bigl(e^{\sqrt{-1}\theta},\,
e^{-\sqrt{-1}\theta}\Bigr):\,\theta\in {\Bbb R}\right\}$$ of ${\rm
SU}(2)$. Hence $G={\tilde G}/\{\pm I_2\}$ is abelian.

Looking at the Riemann sphere $\BCB$ as the complex
one-dimensional projective space ${\Bbb P}(V)$ with the natural
projection $\pi:V\backslash \{0\}\to {\Bbb P}(V)$, we can see that
both $\pi(e_1)$ and $\pi(e_2)$ are two distinct fixed points of
the $G$-action on ${\Bbb P}(V)$ if $G$ can be diagonalized.

Suppose that $G$ has a fixed point $\pi(e_1)$ on ${\Bbb P}(V)$
with $e_1\in V$ and $\langle e_1\,,\,e_1\rangle=1$. Then, $e_1$ is
a common eigenvector of all the elements in ${\widetilde G}$.
Since $\ro_{\widetilde G}$ is a unitary representation on $V$, it
can be decomposed into $V={\Bbb C}e_1\oplus {\Bbb C}e_2$, where
$e_2$ is a unit vector orthogonal to $e_1$. That is, $G$ can be
diagonalized.\\

\nd (2) The abelian subgroup $D_2$ of $\SU$ generated by
\[z\mapsto -z\quad{\rm and}\quad z\mapsto\df{1}{z} \] has no fixed
point on $\BCB$.\hf

\begin{lem}
\label{lem:trivial} A trivial reducible metric $g$ representing a
divisor $D$ on a compact Riemann surface $\Sigma$ is a pullback
$f^*g_{\rm st}$ of $g_{\rm st}$ by some rational function $f$ on
$\Sigma$.
\end{lem}

\nd \pf Let $f:\Sigma^*\to\BCB$ be a developing map of the metric
$g$. Since the monodromy of $f$ is trivial, so is the local
monodromy of $f$ around each point in the support of $D$. By Lemma
\ref{lem:char} the divisor $D$ must be an effective ${\Bbb
Z}$-divisor. Using Lemma \ref{lem:char} again, we find that the
holomorphic map $f:\Sigma^*\to \BCB$ has a holomorphic extension
to $\Sigma$.\hf

\begin{lem}
\label{lem:abelPsi} Suppose that $g$ is an reducible metric on
$\Sigma$ and $f$ a multiplicative developing map of $g$. Then the
holomorphic 1-form $d\,(\log\,f)$ on $\Sigma^*$ can be extended to
be an abelian differential of 3rd kind on $\Sigma$. The function
$\Psi=4|f(z)|^2/(1+|f(z)|^2)$ on $\Sigma^*$ can be extended
continuously to $\Sigma$.
\end{lem}

\nd \pf The proof is contained in the following proof of Theorem
\ref{thm:nontrivial}.\hf\\

\nd {\bf Proof of Theorem \ref{thm:nontrivial}}

\nd {\bf (1-2)}\quad We show that {\it if a point $p\in\Sigma^*$
is a zero of $Y(z)=\frac{f(z)}{f'(z)}\frac{\partial}{\partial z}$,
then $p$ is simple}. We choose a function element ${\frak f}$ near
$p$. Since the monodromy of $f$ belongs to ${\rm U}(1)$,
$Y=\frac{{\frak f}(z)}{{\frak f}'(z)}\frac{\partial}{\partial z}$,
which is independent of the choice of the function element ${\frak
f}$ and the complex coordinate $z$. Since ${\frak f}$ is a
univalent meromorphic function near $p$, there exists ${\frak
L}\in \PGL$ and a complex coordinate $z$ near $p$ with $z(p)=0$
such that ${\frak L}\circ f=z$. Then ${\frak f}=\frac{az+b}{cz+d}$
with $ad-bc=1$ near $p$, and
$Y=(az+b)(cz+d)\frac{\partial}{\partial z}$. It is clear that $p$
could not be a pole of $Y$. Since $Y=0$ at $z(p)=0$, $bd=0$.

\nd {\it Case 1}\quad Since $ad-bc=1$ and $bd=0$, we assume $b=0$,
$d\not=0$ in this case. Then $ad=1$, ${\frak
f}(z)=\frac{az}{cz+d}$. Then $Y=az(cz+d)\frac{\partial}{\partial
z}$ has a simple zero at $z(p)=0$. Hence
$\omega=\frac{dz}{az(cz+d)}$ has residue $1$ at $p$. Since
$f(0)=0$, $p$ is a minimal point of $\Psi$ and $\Psi(p)=0$.

\nd {\it Case 2}\quad Similarly, when $d=0$ and $b\not=0$, $Y$
also has simple zero at $p$, $\omega$ has residue $-1$,
$\lim_{q\to p}\, |f(q)|=+\infty$,  $\lim_{q\to p}\,\Psi(q)=4$, and
$p$ is a maximal point of $\Psi$.  \\

We show that {\it each point $q\in \{p_1,\cdots,p_n\}$ must be a
simple zero of $Y$, provided the conical angle of the metric $g$
at $q$ equals $2\pi s>0$ and $s$ is a non-integer}. By Lemmata
\ref{lem:Sch} and \ref{lem:char}, we can choose a function element
${\frak f}$ near $q$ and a complex coordinate $z$ near $q$ such
that ${\frak f}=\frac{az^s+b}{cz^s+d}$ with $ad-bc=1$. On the
other hand, since the monodromy of $f$ belongs to ${\rm U}(1)$, so
does the local monodromy of ${\frak f}$. Then there exists
$\theta\in{\Bbb R}$ such that
\[e^{2\pi\sqrt{-1}\theta}\,{\frak f}=e^{2\pi\sqrt{-1}\theta}\,\frac{az^s+b}{cz^s+d}=
\df{a\, e^{2\pi\sqrt{-1} s}\,z^s+b}{c\, e^{2\pi\sqrt{-1}
s}\,z^s+d}. \] This is equivalent to that the following equalities
hold,
\begin{eqnarray*}
ac\,e^{2\pi\sqrt{-1}s}\left(1-e^{2\pi\sqrt{-1}\theta}\right)&=&0,\\
\left(ad\,e^{2\pi\sqrt{-1}s}+bc\right)-
e^{2\pi\sqrt{-1}\theta}\,\left(bc\, e^{2\pi\sqrt{-1}s}+ad\right)
&=&0,\\
bd\,\left(1-e^{2\pi\sqrt{-1}\theta}\right)&=&0.
\end{eqnarray*}
Solving the equation, we find that either $c=b=0$ or $a=d=0$. i.e.
${\frak f}(z)$ equals $\mu\, z^s$ ($\mu\not=0$) or $\lambda\,
z^{-s}$ ($\lambda\not=0$). Hence $Y=\pm
s\,z\,\frac{\partial}{\partial z}$ has simple zero at $z(q)=0$.
Since ${\frak f}(0)$ equals $0$ or $\infty$, $\Psi$ is continuous
at $q$, which is a minimal or maximal point of $\Psi$ achieving
value $0$ or $4$ if and only if $\omega$
has residue $s>0$ or $-s<0$ at $p$.\\

Let $p$ be a singular point of the metric $g$ with conical angle
$2\pi$ times an integer $n>1$. Then ${\frak
f}(z)=\df{az^n+b}{cz^n+d}$ and
$Y=\df{(az^n+b)(cz^n+d)}{nz^{n-1}}\,\df{\partial}{\partial z}$.

\nd{\it Case A}\quad Assume $bd\not=0$. Then $p$ is a pole of $Y$
with order $n-1$ and a zero of $\omega$ with order $n-1$, and
$\lim_{z\to p}\,f(z)=f(0)=b/d\in {\Bbb C}\backslash\{0\}$.
Moreover, $\Psi$ is continuous at $p$, which is a saddle point of
$\Psi$.

\nd {\it Case B}\quad Assume $bd=0$. Then it is easy to check that
$p$ is a simple zero of $Y$. If $b=0$ ($b\not=0$), then
$\lim_{q\to p}\,|f(p)|$ equals $0$ or $+\infty$, where $\omega$
has residue $n$ or $-n$. $\Psi$ is continuous at $p$, where it
achieves the minimal value $0$ or the maximal  value $4$.\\

\nd {\bf (3)}\quad The local monodromy property of $f$ follows
from
Lemmata \ref{lem:Sch} and \ref{lem:char}. \hf\\

\begin{lem}
\label{lem:pole} The following two statements hold.

\nd {\rm (1)} Let $g$ be an reducible metric on $\Sigma$. Then
each character 1-form $\omega$ of $g$ has at least two poles.

\nd {\rm (2)} Besides the assumption of (1), assume that $\omega$
has no zero and has only two poles. Then $\Sigma$ is the Riemann
sphere $\BCB$ and $g$ has two singularities with the same angle,
say $\alpha>0$. Moreover, if the two singularities are assumed to
be $0$ and $\infty$, then $\omega=\alpha\,\df{dz}{z}$ up to sign.
\end{lem}

\nd\pf (1) Let $f$ be the multiplicative developing map such that
$\omega=d\,(\log\,f)$. Since $\Psi=4|f|^2/(1+|f|^2)$ is a
non-constant continuous function on $\Sigma$, it must achieve
minimum and maximum. Either a minimal point or a maximal one of
$\Psi$ is a pole of $\omega$ by Theorem \ref{thm:nontrivial}.

(2) $\Sigma=\BCB$ follows from $\deg\,(\omega)=-2$. By the Residue
Theorem, the two residues of $\omega$ have the different sign at
the two poles, say $0$ and $\infty$. It follows from Theorem
\ref{thm:nontrivial} that $g$ has exactly two singularities $0$
and $\infty$ with the same angle, say $\alpha$. Assume that
$\omega$ has residue $\alpha$ ($-\alpha$) at $0$ ($\infty$). Then
$\omega=\alpha\,\df{dz}{z}$.
\hf\\

\nd {\bf Proof of Theorem \ref{thm:char}}\quad We divide the proof
into the
following two cases.\\

{\it Case 1} \quad Assume that the integral of $\omega$ at some
loop in $\Sigma':=\Sigma\backslash\{\text{poles of }\omega\}$ does
not belong to the set $2\pi\sqrt{-1}\,{\Bbb Z}$. Since
$\Re\,\omega$ is exact on $\Sigma'$, sloving the equation
$$\omega=d\,(\log\,f)$$ on $\Sigma'$, we obtain a multi-valued
locally univalent meromorphic function
\[f(z)=\exp\,\left( \int^z\, \omega \right)\]
unique up to a complex multiple with modulus one. Moreover, $f$
has non-trivial monodromy belonging to ${\rm U}(1)$ and
$f^*\,g_{\rm st}$ is an non-trivial reducible metric with
character 1-form $\omega$. Conversely, if $g$ is an reducible
metric such that $\omega$ is one of its character 1-forms, then
there exists a developing map ${\tilde f}$ of $g$ such that
$\omega=d\,(\log\, {\tilde f})$. Since ${\tilde f}$ has
non-trivial monodromy in ${\rm U}(1)$, $g$ is non-trivial. Solving
the equation $\omega=\frac{d\,{\tilde f}}{\tilde f}$ also gives
the same expression of ${\tilde f}$ as $f$. Therefore, such
reducible metric $g$ is unique. By the argument in the proof of
Theorem \ref{thm:nontrivial}, we find the divisor $D$ represented
by $g$ equals
\[\sum_{j=1}^J\,(\alpha_j-1)\,P_j+
\sum_{k=J+1}^N\,\Bigl(|{\rm Res}_{Q_k}(\omega)|-1\Bigr)\, Q_k.\]

{\it Case 2}\quad Assume that the monodromy given by $\omega$ is
trivial, i.e. the integral of $\omega$ at each loop in
$\Sigma':=\Sigma\backslash\{\text{poles of }\omega\}$ belongs to
the set $2\pi\sqrt{-1}\,{\Bbb Z}$. $f^*\,{\rm g}_{\rm st}$ with
$f(z)=\exp\,\left( \int^z\, \omega \right)$ is a trivial reducible
metric such that $f$ is one of its developing map and a rational
function on $\Sigma$ and $\omega$ is one of its character 1-forms.
Conversely, if $g$ is an reducible metric with $\omega$ one of its
character 1-forms, then $g={\tilde f}^*\,g_{\rm st}$ with ${\tilde
f}(z)=\exp\,\left( \int^z\, \omega \right)$. Moreover, ${\tilde
f}$ is a rational function uniquely determined by $\omega$.
Therefore, such reducible metric $g$ is unique. By the similar
argument in the proof of Theorem \ref{thm:nontrivial}, we can show
that the effective divisor represented by $g$ equals
$\sum_{j=1}^J\,(\alpha_j-1)\,P_j+ \sum_{k=J+1}^N\,\Bigl(|{\rm
Res}_{Q_k}(\omega)|-1\Bigr)\, Q_k$.
\hf\\

\begin{cor}
\label{cor:tro} Under the notations in Theorem \ref{thm:char}, we
have
\[\chi(\Sigma)+\deg\, D\geq \min\,\bigl(2,\,2\,\min_j\,\alpha_j\bigr).\]
In particular, the divisor $D$ does not satisfy Troyanov's
condition \eqref{equ:tro}. In other words, if $D$ satisfies
condition \eqref{equ:tro}, then each conformal metric, which has
constant curvature one and represents $D$, is irreducible.
\end{cor}

\nd\pf The character 1-form $\omega$ has at least two poles since
the continuous function $\Psi=4|f|^2/(1+|f|^2)$ in Theorem
\ref{thm:nontrivial} has at least a minimal point and a maximal
one.  Then, using the equality
\[-\chi(\Sigma)=\deg\,(\omega)=\sum_{j=1}^J\,(\alpha_j-1)-(N-J),\]
we have
\begin{eqnarray*}
\chi(\Sigma)+\deg\, D&=&\chi(\Sigma)+\sum_{j=1}^J\,
(\alpha_j-1)+\sum_{k=J+1}^N\,\bigl(|{\rm
Res}_{Q_k}\,(\omega)|-1\bigr)\\
&=&\sum_{k=J+1}^N\,|{\rm Res}_{Q_k}\,(\omega)|\geq
\min\,\bigl(2,\,2\min_j\,\alpha_j\bigr).
\end{eqnarray*}
\hf\\

\nd \begin{exam} \label{exam:nonabel}

Consider on the two-sphere a conformal metric $g$ with constant
curvature one and finite conical singularities
$p_1,\,\cdots,\,p_n$. Let the angle at $p_j$ be $2\pi\,\alpha_j$.
If $n\geq 3$ and each $\alpha_j$ is a non-integer, then $g$ is
irreducible. Otherwise, by Theorem \ref{thm:nontrivial}, the
character 1-form of $g$ would have at least three poles and have
no zeroes. Contradiction. In particular, we consider the conformal
metric $g$ of constant curvature one with three angles
$\pi,\,\pi,\,\pi$ at $0,\,1$ and $\infty$ on the two-sphere. Then
the developing map of $g$ is the Gauss hypergeometric function $u$
satisfying the hypergeometric equation
\[z(1-z)\df{d^2
u}{dz^2}+\bigl(c-(a+b+1)z\bigr)\df{du}{dz}-abu=0, \quad {\rm
where}\quad |1-c|=|c-a-b|=|a-b|=\df{1}{2}.\] Then the monodromy
group of a developing map of $g$ is conjugate to $D_2$ in Lemma
\ref{lem:diag} (cf \cite[p. 53]{Yo87}) and then abelian, but $g$
is irreducible.

It follows from Troyanov \cite[Theorem 4]{Tr91} that on a torus
there exists a conformal metric $g$ with constant curvature one
and a conical singularity $p$ with angle $2\pi\,\alpha$, where
$1<\alpha<3$. Then Corollary \ref{cor:tro} tell us that  $g$ is
not reducible. The existence of such a irreducible metric having
one angle $4\,\pi$ implies that Theorem \ref{thm:rational} does
not hold on a torus.
\end{exam}

\begin{exam}
\label{exam:abel} Let $a,\,b$ be two positive numbers. Consider
the 1-form
\[\omega=\left(\df{a}{z}+\df{b}{z-1}\right)\,dz,\]
which has residue $-a-b$ at $\infty$. $a/(a+b)$ is the zero of
$\omega$. Hence $\omega$ satisfies the condition in Theorem
\ref{thm:char}, and then gives an reducible metric on the
two-sphere with angles $2\pi a,\,2\pi b,\,2\pi (a+b)$ and $4\pi$
at $0,\,1,\,\infty$ and $a/(a+b)$, respectively.

Suppose that $g$ is an reducible metric on the two-sphere having
angles $2\pi a,\,2\pi b,\,2\pi (a+b)$ and $4\pi$ at
$0,\,1,\,\infty$ and $\lambda\in {\Bbb C}\backslash\{0,\,1\}$,
respectively, where none of $a,\,b,\, a+b$ is an integer. Then
$\lambda=a/(a+b)$. Actually, letting $\omega$ be a character
1-form of $g$, we can see from Theorem \ref{thm:nontrivial} that
$0,\,1,\,\infty$ are simple poles of $\omega$ and $\lambda$ is the
zero of $\omega$. By the Residue Theorem, we may assume that
\[{\rm Res}_{0}(\omega)=a,\quad {\rm Res}_{1}(\omega)=b,\quad
{\rm Res}_{\infty}(\omega)=-a-b,\] which implies that
$\omega=\left(\df{a}{z}+\df{b}{z-1}\right)\,dz$ and
$\lambda=\df{a}{a+b}$. This implies that the existence of
reducible metrics does not only depend on angles, but also on the
position of singularities.
\end{exam}

\nd {\bf Remark 4.1.} The monodromy of a developing map is
irreducible for a hyperbolic conformal metric with finite conical
or cusp singularities. Moreover, it is also the case for a flat
conformal metric with finite conical singularities unless the
metric is a smooth one on a torus. We leave the proof to the
interested readers.\\

We conclude this section by saying something more about the
relationship between reducible metrics and HCMU metrics. As a
generalization with singularities on compact Riemann surfaces of
Calabi's extremal K{\" a}hler metric on compact complex manifolds
(\cite{Ca82, Ca85}), X. Chen \cite{Chen00, Chen98, Chen00} firstly
introduced the concept of HCMU metric and extremal metric with
singularities and proved some fundamental results. In particular,
a conformal metric ${\tilde g}$ on a compact Riemann surface with
singularities is called {\it HCMU} if and only if it has finite
area and finite Calabi energy (cf Section 5), and the complex
gradient $K^{,\,z}\,\tpds{}{z}$ of the Gauss curvature
$K=K_{\tilde g}$ is a holomorphic vector field outside the
singularities. Wang-Zhu \cite{WZ00} and Lin-Zhu \cite{LZ02}
obtained some interesting results and generalized some results of
X. Chen. Recently,  three of the authors \cite{CW11, CWX12} made
the complete classification of non-constant curvature HCMU metrics
with conical or cusp singuarities, by using the {\it character
1-form}
$${\tilde
\omega}=\df{dz}{K_{\tilde g}^{,\,z}}$$ of a HCMU metric ${\tilde
g}$. The property plays a crucial rule in the classification that
${\tilde \omega}$ is an abelian differential of 3rd kind with real
residues and its real part is exact outside the set of simple
poles. The following observation we made lay both HCMU metrics and
metrics of constant curvature
in the same philosophical frame. \\

\nd {\bf Observation} {\it Given a non-constant curvature HCMU
metric ${\tilde g}$ with singularities on a compact Riemann
surface $\Sigma$, there exists a multi-valued locally univalent
meromorphic function ${\tilde f}$ on $\Sigma^*$ having monodromy
in the abelian group
$$\left\{\exp\,\bigl(\sqrt{-1}\theta\bigr)|\theta\in {\Bbb
R}\right\},$$ and a football HCMU metric $g_{\rm fb}$ (cf
\cite{Chen00, CCW05}) over $\BCB$ such that ${\tilde g}={\tilde
f}^*g_{fb}$. Moreover, the character 1-form of ${\tilde g}$
coincides with the logarithmic differential
$d\,(\log\, f)$ of $f$, up to a constant.}\\

\nd It is not used in this paper and its proof will be left
elsewhere.

\section{Discussions}
\label{sec:appl}

As an application of Theorems \ref{thm:rational} and
\ref{thm:nontrivial}, we shall show that {\it if $g$ is conformal
metric on the sphere $\BCB$ of constant curvature one and
representing the divisor $D=(\alpha-1)P+(\beta-1)Q$, where
$\alpha,\,\beta>0$, then $\alpha=\beta$.}

\nd \pf {\it Case 1}\quad We assume that at least one of $\alpha$
and $\beta$ is not an integer. Suppose that it is the case for
$\alpha$. Since the punctured sphere $\BCB\backslash\{p,\,q\}$ has
the fundamental group isomorphic to ${\Bbb Z}$, the metric $g$ is
a reducible metric. Let $f$ be one of its developing map. By
Lemmata \ref{lem:Sch} and \ref{lem:char}, the local monodromy of
$f$ at $p$ is non-trivial. Hence $g$ is non-trivial.  We may
assume that $f$ is multiplicative so that $\omega=\frac{df}{f}$ is
the character 1-form of $g$. Theorem \ref{thm:nontrivial} tells us
that $p$ is a simple pole of $\omega$ with residue $\pm \alpha$
and $q$ is either a simple pole or a zero point of $\omega$.  If
$q$ is a simple pole too, then the residue equals $\pm \beta$.
Since the only poles of $\omega$ are $p,\,q$, we have
$\alpha=\beta$. We shall rule out the case where $q$ is a zero
point of $\omega$. Suppose it is the case. Then, $\omega$ has only
one simple pole in $\Sigma^*$, which has residue $\pm 1$ by
Theorem \ref{thm:nontrivial}. It contradicts the fact that
$\alpha\not=1$.

{\it Case 2}\quad Suppose that both $\alpha$ and $\beta$ are
integers $\geq 2$. By Theorem \ref{thm:rational}, each developing
map $f$ of $g$ is a rational function on $\BCB$. We may assume
that $p=0$ and $f(p)=0$ by using suitable fractional linear
transformations. Since $f$ is rational, $f(z)=z^\alpha$ in some
complex coordinate $z$ near $0$, which is a simple pole of
$\omega=\frac{df}{f}=\frac{\alpha\,dz}{z}$. Similarly, $q$ is also
a simple pole of $\omega$. The residue theorem gives $\alpha=\beta$.\hf \\

We observe that the reducible metrics representing a given divisor
are not unique in general in the following

\begin{prop}
\label{prop:nonuni} Suppose that there exists an reducible metric
$g$ representing the divisor $D$ on the two-sphere. Denote by
${\Bbb M}(D)$ the space of conformal metrics of constant curvature
one representing $D$, by ${\Bbb A}(D)$ that of reducible
metrics representing $D$. \\

\nd {\bf (1)} If $D$ is supported at two points $p_1$ and $p_2$,
${\Bbb M}(D)={\Bbb A}(D)$, and $g$ is unique if and only if the
two angles are equal and do not belong to $2\pi\,{\Bbb Z}_{>1}$.
If the two angles are equal and belong to $2\pi\,{\Bbb Z}_{>1}$,
then ${\Bbb A}(D)$ is connected and has dimension
1.\\

\nd {\bf (2)} If $D$ is supported at three points and ${\Bbb
A}(D)\not=\emptyset$, then ${\Bbb M}(D)={\Bbb A}(D)$ is connected.
Moreover, if $g$ is trivial,
$\dim\,{\Bbb A}(D)=3$; otherwise $\dim\,{\Bbb A}(D)=1$.\\

\nd {\bf (3)} Suppose that $D$ is supported at more than three
points and ${\Bbb A}(D)\not=\emptyset$. Then, if $g$ is trivial,
$\dim\,{\Bbb A}(D)=3$; otherwise $\dim\,{\Bbb A}(D)\geq 1$.\\

\end{prop}

\nd \pf The first statement was proved by Troyanov \cite[Theorem
I, p. 298]{Tr89}. The second was shown in Umehara-Yamada
\cite[Corollary 2.3]{UY00}.

Suppose that $D$ supports at more than three points. Following
Umehara-Yamada \cite[(2.5)]{UY00}, we define
\[I_g:=\left\{g_a=(a\star f)^*\,g|a\in {\rm PSL}(2,\,\BC);\, a\cdot
{\rm Im}\,\ro_f\cdot a^{-1}\subset \SU\right\},\] where $a\star f$
means the M{\" o}bius transformation of $f$ by $a$, and
$$\ro_f:\pi_1(\Sigma^*)\to \SU$$ denotes the monodromy
representation of the developing map $f$ of the metric $g$. Each
metric $g_a$ in $I_g$ has a developing map $a\star f$, which has
the same Schwarzian with $f$ and monodromy conjugate to that of
$f$. Hence $I_g$ is contained in ${\Bbb A}(D)$. Then it follows
from \cite[Lemma B, p. 92]{UY00} that if $g$ is trivial,
$\dim\,{\Bbb A}(D)\geq 3$; otherwise $\dim\,{\Bbb A}(D)\geq 1$.

We consider the moduli of trivial reducible metrics representing
an effective ${\Bbb Z}$-divisor $D$, which can be reduced to the
space of rational functions with the same ramification divisor
$D$. We say that two rational functions have the {\it same type}
if one of them is given by the post-composition with a M{\" o}bius
transformation of the other. It follows from \cite[Lemma B, p.
92]{UY00} that the trivial reducible metrics having developing
maps of the same type form a moduli of the three dimensional
hyperbolic space ${\cal H}^3$. The beautiful theorem in I.
Scherbak \cite{Sch02} says that there is a least upper bound given
by the Schubert calculus for the number of types of all the
rational functions with ramification divisor $D$, which can be
achieved by a generic choice of the support of $D$. Hence we
obtain the corresponding information for the number of connected
components of ${\Bbb
A}(D)$. \hf\\

It is time for us to propose some questions interesting to us. \\

\nd {\bf Question 1.} Does there exist a divisor on some compact
Riemann surface, which is represented by both an irreducible
metric and a reducible one?  It does not happen on the two-sphere
under either of the following two conditions:

\nd (i) The support of $D$ consists of three points or less (cf
Proposition \ref{prop:nonuni}).

\nd (ii) Each $\alpha_j$ is a non-integer (cf Example
\ref{exam:nonabel}).\\

\nd {\bf Question 2.} Suppose ${\Bbb A}(D)$ is non-empty for a
divisor $D$ on a compact Riemann surface. Study the moduli space
${\Bbb A}(D)$ of reducible metrics representing a divisor $D$ on a
compact Riemann surface, such as its dimension and the number of
its components. We know the answer on the two-sphere in case that
$D$ supports at two or three points or $D$ is an effective ${\Bbb
Z}$-divisor in Proposition
\ref{prop:nonuni}. \\
\\

\nd {\bf Question 3.} Suppose that there exists an irreducible
metric $g$ representing $D=\sum_j\,(\alpha_j-1)\,P_j$. Is $g$ the
unique metric of constant curvature one representing $D$? Luo-Tian
\cite{LT92} showed it is the case on the two-sphere if each
$\alpha_j$ lies in $(0,\,1)$. Moreover  Umehara-Yamada \cite{UY00}
gave the positive answer if $D$ is a divisor supporting at three
points on the two-sphere.\\
\\

\nd{\bf Question 4.} Suppose that there exists an irreducible
metric $g$ representing $D=\sum_j\,(\alpha_j-1)\,P_j$. Does there
exist an irreducible metric representing any divisor $D'$
sufficiently near $D$?  On the two-sphere, if each $\alpha_j$ lies
in $(0,\,1)$, then the necessary and sufficient condition is an
open one for the existence of a irreducible metric representing
$D$ given by Troyanov \cite{Tr91} and Luo-Tian \cite{LT92}. On the
two-sphere, if $D$ supports at three points, so is the necessary
and sufficient condition given by Umehara-Yamada \cite{UY00}. S.
K. Donaldson \cite{Don11} proved an openness theorem for K{\"
a}hler Einstein metrics on a Fano manifold with conical
singularity along the anti-canonical divisor.

\vspace{0.5cm}

\section{Cusp singularity}

Let $g=e^{2\varphi}\,|dz|^2$ be a smooth conformal metric in a
punctured disk $D\backslash\{0\}$ with finite area and finite
Calabi energy, where the Calabi energy $E(g)$ of $g$ is defined by
the square integral of the Gauss curvature over
$(D\backslash\{0\},\,g)$. Then X. Chen proved in \cite[Theorem 2,
p.198]{Chen98} that
\begin{equation}
\label{equ:Chen} \lim_{\ro\to 0}\,
\left(\varphi(re^{i\theta})+\ln\,r\right)=-\infty,
\end{equation}
where $z=re^{i\theta},\,r=|z|$. We call that $0$ is a {\it genuine
weak cusp singular point} (cf \cite[p.195]{Chen98})of the
conformal metric $g$ if $\varphi$ satisfies the following integral
condition
\begin{equation}
\label{equ:wcusp} \liminf_{r\to 0}\,
\int_0^{2\pi}\,r\,\tpds{(\varphi+\ln\,r)}{r}=0.
\end{equation}

\begin{thm}
\label{thm:cusp} {\rm (\cite[Theorem 2]{Wang11})}  If the
aforementioned metric $g=e^{2\varphi}\,|dz|^2$ has nonnegative
Gauss curvature in the punctured disk $D\backslash\{0\}$, then $0$
could not be its genuine weak cusp singular point.
\end{thm}
 \nd
\pf The proof of the theorem is essentially the same as Lemma 6 in
\cite[p.218]{Chen99}. Under the coordinate transformation $t=\ln\,
r,\,\theta=\theta$, the punctured disk $D\backslash
\{0\}=\{z=r\,e^{i\theta}:\,0<r<1,\,-\pi\leq\theta\leq\pi\}$ is
transformed to the infinite cylinder
$$\{(t,\,\theta):\, -\infty<t<0,-\pi\leq\theta\leq\pi\}.$$
Then the metric $g=e^{2\varphi}\,|dz|^2$ has the new expression
$$g=e^{2\psi}\,|dw|^2,\quad
\psi(t+i\theta)=\varphi(re^{i\theta})+t,$$ under the new complex
coordinate $w=t+i\theta$. \eqref{equ:Chen} means that
$\lim_{t\to-\infty}\,\psi(t+i\theta)=-\infty$ for each $\theta\in
[-\pi,\,\pi]$. Then
$\Psi(t):=\int_0^{2\pi}\,\psi(t+i\theta)\,d\theta$ diverges to
$-\infty$ as $u$ goes to $-\infty$. Then we could choose $t_1<0$
with $\Psi'(t_1)>0$, i.e.
\begin{equation}
\label{equ:Chen2} \int_0^{2\pi}\,
\tpds{\psi}{t}(t_1+i\theta)\,d\theta>0.
\end{equation}
Since the metric $g=e^{2\psi}\,|dw|^2$ has nonnegative Gauss
curvature $K$ in the cylinder, we have
\[-\Delta_{t,\,\theta}\,\psi(t+i\theta)=K\,e^{2\psi}\geq 0\quad
\text{for all}\quad t<t_1.\] By the Green formula, integrating the
above equation in $[t,\,t_1]\times S^1$, we obtain
\[0\leq -\int_{[t,\,t_1]\times S^1}\,
\Delta_{t,\,\theta}\,\psi(t+i\theta)\,dt\,d\theta
=\Psi'(t)-\Psi'(t_1).\] By \eqref{equ:Chen2},
\begin{eqnarray*}
\liminf_{r\to 0}\,
\int_0^{2\pi}\,r\,\tpds{(\varphi+\ln\,r)}{r}(re^{i\theta})\,d\theta
&=&\liminf_{t\to-\infty}\,\int_0^{2\pi}\,
\tpds{\psi}{t}(t+i\theta)\,d\theta\\
&=&\liminf_{t\to-\infty}\,\Psi'(t) \geq \Psi'(t_1)>0.
\end{eqnarray*}
Therefore $0$ could not be the genuine weak cusp singularity of
the metric $g$ by definition. \hf\\

\nd{\bf Remark} We call that $0$ is a {\it cusp singularity} of
the metric $g=e^{2\varphi}\,|dz|^2$ if there holds
\[\lim_{r\to\infty}\,\frac{\varphi+\ln\,r}{\ln\,r}=0.\]
It is proved in \cite{CWX12} that $0$ is a genuine weak cusp
singularity of $g$ if and only if it is a cusp singularity.
Intuitively, we may think of a cusp singularity as a conical one
with angle $0$.\\

%\nd ${\frak A}\ {\frak B}\ {\frak C}\ {\frak D}\ {\frak E}\ {\frak
%F}\ {\frak G}\\ {\frak H}\ {\frak I}\ {\frak J}\ {\frak K}\ {\frak
%L}\ {\frak M}\ {\frak N}\\ {\frak O}\ {\frak P}\ {\frak Q}\ {\frak
%R}\ {\frak S}\ {\frak T}\ {\frak U}\ {\frak V}\ {\frak W}\ {\frak
%X}\ {\frak Y}\ {\frak Z}$
%\\

%\nd ${\frak a}\ {\frak b}\ {\frak c}\ {\frak d}\ {\frak e}\ {\frak
%f}\ {\frak g}\\ {\frak h}\ {\frak i}\ {\frak j}\ {\frak k}\ {\frak
%l}\ {\frak m}\ {\frak n}\\ {\frak o}\ {\frak p}\ {\frak q}\ {\frak
%r}\ {\frak s}\ {\frak t}\ {\frak u}\ {\frak v}\ {\frak w}\ {\frak
%x}\ {\frak y}\ {\frak z}$

\end{document}